\documentclass[review]{elsarticle}

\usepackage{lineno,hyperref}
\modulolinenumbers[5]

\journal{}

\usepackage{amsmath}
\usepackage{amsfonts}
\usepackage{amssymb}

\usepackage{graphicx}
\usepackage{latexsym}
\usepackage{amsthm}
\usepackage{xcolor}

\newcommand{\R}{\mathbb{R}} %Real numbers

\newcommand{\be}{\begin{equation}}
\newcommand{\ee}{\end{equation}}
\newtheorem{theorem}{Theorem}[section]

\def\Proof{{\noindent \bf Proof.\enspace}}
\def\Endproof{\indent\vrule height6pt width4pt depth1pt\hfil\par\medbreak}
\begin{document}
\begin{frontmatter}

\title
 {Modified Singly-Runge-Kutta-TASE
 methods for the numerical solution of stiff differential equations\tnoteref{mytitlenote}}

 \tnotetext[mytitlenote]{Partially supported by project PID2022-141385NB-I00 of Ministerio de Ciencia e Innovaci\'{o}n of Spain.}

\author[Zamainaddress]{M. Calvo}
\ead{calvo@unizar.es}

\author[Zamainaddress]{J. I. Montijano\corref{mycorrespondingauthor}}
\cortext[mycorrespondingauthor]{Corresponding author}
\ead{monti@unizar.es}

\author[Zamainaddress]{L. R\'{a}ndez}
\ead{randez@unizar.es}

\address[Zamainaddress]{Departamento  Matem\'{a}tica Aplicada,
Universidad de Zaragoza.
50009-Zaragoza, Spain.}

\begin{abstract}
Singly-TASE operators for the numerical solution of stiff differential equations were proposed by Calvo et al. in J.Sci. Comput. 2023 to reduce the computational cost of Runge-Kutta-TASE (RKTASE) methods when the involved linear systems are solved by some $LU$ factorization. In this paper we propose a modification of these methods to improve the efficiency by considering  different TASE operators for each stage of the Runge-Kutta. We prove that the resulting RKTASE methods are equivalent to $W$-methods (Steihaug and Wolfbrandt, Mathematics of Computation,1979) and this allows us to obtain the order conditions of the proposed Modified Singly-RKTASE
methods (MSRKTASE) through the theory developed for the $W$-methods. We construct new MSRKTASE methods of order two and three and demonstrate their effectiveness through numerical experiments on both linear and nonlinear stiff systems. The results show that the  MSRKTASE schemes significantly enhance efficiency and accuracy compared to previous Singly-RKTASE schemes.
\end{abstract}


\end{frontmatter}

\section{ Introduction}

Solving stiff initial value problems (IVPs) efficiently and accurately remains a significant challenge in numerical analysis. Explicit Runge-Kutta (RK) method, face limitations when applied to stiff systems due to their stability restrictions.
To address these limitations, a new class of time-advancing schemes for the numerical solution of stiff IVPs was proposed by M. Bassenne, L. Fu, and A. Mani in \cite{Bassenne:2020}. They introduced the concept of TASE (Time-Accurate and Stable Explicit) operators, designed to enhance the stability of explicit RK methods. In this approach, instead of solving
the differential system
\begin{equation}\label{A1}
\dfrac{d}{dt} Y(t) = F ( t, Y(t)),
\quad Y(t_0)=Y_0 \in \R^d.
\end{equation}
they proposed to solve another (stabilized) IVP
\begin{equation}\label{A2}
\dfrac{d}{dt} U(t) = T \;  F ( t, U(t)),
\quad U(t_0)=U_0 \equiv Y_0 \in \R^d,
\end{equation}
 where $ T = T( \Delta t)$ is a linear operator that may depend
on the time step size $\Delta t >0$  and on $F$,
so that the numerical solution of Eq.~\eqref{A2} $U_{RK}(t_0 + \Delta t)$
obtained with an explicit RK method of order $p$ satisfies
$
 U_{RK}(t_0 + \Delta t) - Y(t_0 + \Delta t) = \mathcal{O}(\Delta t ^{p+1}),
$
i.e., approximates the local solution of Eq.~\eqref{A1} with order $p$, and also
satisfies some stability requirements that are necessary for solving stiff systems
such as  A- or L- stability.
In this way, the introduction of the TASE operator into the
original governing equation \eqref{A1} allows to overcome the numerical stability restrictions of explicit RK time-advancing methods for solving stiff systems.

In terms of the accuracy order of the numerical solution of Eq.~\eqref{A2}, if the explicit RK scheme has order $p$
and if the TASE operator satisfies
\begin{equation}\label{A3}
  T ( \Delta t ) = I + \mathcal{O} ( \Delta t^q),
\end{equation}
it can be easily proved that the numerical solution of the stabilized system Eq.~\eqref{A2}, $U_{RK}(t_n)$, has at least order $\min(p,q)$.
Some specific schemes combining explicit $p$--stage RK methods with orders $ p \le 4$ and TASE operators with the same order were derived in \cite{Bassenne:2020}.

A more general family of TASE operators was derived in \cite{CaMoRa:21} by taking
\begin{equation}\label{A8}
T(\Delta t)=
\sum_{j=1}^p \beta_j (I-\alpha_j\Delta t W)^{-1} ,
\end{equation}
where $ \alpha_j >0$, $ j=1, \ldots ,p$, are free real parameters, and $\beta_j$ are uniquely determined by the condition
$ T_p ( \Delta t ) = I + \mathcal{O} ( \Delta t^p) $.
The free parameters $\alpha_j$ were selected to improve the linear stability properties of an explicit
RK method with order $p$ for the scaled equation Eq.~\eqref{A2}.

In the TASE operators described by Eq.~(\ref{A8}), each evaluation of $T$
involves the solution of $p$ linear systems with matrices
$ ( I - \Delta t \; \alpha_j \; W)$, $ j=1, \ldots ,p$. If these systems
are to be solved by some $LU$ matrix factorization, it would be much more efficient if the $\alpha_j$ coefficients were all equal, but this is not compatible with order $p$ for the TASE operator.

To make the methods more efficient, Calvo et al. \cite{CaFuMoRa:23} developed a different family of TASE operators, termed Singly-TASE operators, defined by
\begin{equation}\label{A9}
T_p(\Delta t)=
\sum_{j=1}^p \beta_j (I-\alpha \Delta t W)^{-j}.
\end{equation}
These operators lead to significant improvements in computational efficiency while maintaining the stability and accuracy necessary for solving stiff problems.

For an $s$--stage explicit RK method defined by the Butcher tableau
\[
\begin{tabular}{c|ccccc}
0 &  & & & &
\\
$c_2 $ & $a_{21} $& & & &
\\
\vdots & & & & &
\\
$c_s$  & $a_{s1}$ & $ a_{s2}$ & \ldots & $ a_{s,s-1}$ &
\\
\hline
 & $b_1$ & $b_2$ & \ldots & $ b_{s-1}$ & $ b_s$
 \end{tabular}
 \qquad
 c_i = \sum_{j=1}^{i-1} a_{ij} , \quad i= 2, \ldots, s,
 \]
the numerical solution of Eq.~(\ref{A2}) by a Singly-RKTASE method  is advanced from  $ (t_0, U_0) \to (t_1= t_0 + \Delta t, U_1)$   by the formula
\[
U_1 = U_0 + \Delta t \left[
b_1 \; K_1 + \ldots + b_s \; K_s \right],
\]
where $ K_j $, $ j=1, \ldots ,s$   are computed recursively from the formulas
\begin{eqnarray}
K_1 &=& T \; F(t_0, U_0),
\nonumber
\\
K_2 &=& T \; F(t_0 + c_2 \Delta t, U_0 + \Delta t \; a_{21} K_1),
\nonumber
\\
\vdots &&
\nonumber
\\
K_s &=&
T \; F \left(t_0 + c_s \Delta t, U_0 + \Delta t \; \sum_{j=1}^{s-1} a_{s j} K_j  \right).
\end{eqnarray}

In this paper, we propose a further modification of the Singly-RKTASE methods to enhance efficiency, so that for the ith stage of the RK method we use a TASE operator defined by
\begin{equation}\label{A10}
T_i = \sum_{j=1}^r \beta_{ij} \; ( I - \Delta t \; \alpha \; W)^{-j},
\quad
\alpha >0.
\end{equation}
The coefficients $\alpha>0$, $r$ and $ \beta_{i,j}$, $ j=1, \ldots ,r$ must be determined so that
the resulting modified singly-RKTASE (MSRKTASE) method
\begin{equation}\label{A11}
\begin{array}{lll}
K_1 &=& T_1 \; F(t_0, U_0),
\\
K_2 &=& T_2 \; F(t_0 + c_2 \Delta t, U_0 + \Delta t \; a_{21} K_1),
\\
\vdots &&
\\
K_s &=&
T_s \; F \left(t_0 + c_s \Delta t, U_0 + \Delta t \; \sum_{j=1}^{s-1} a_{s j} K_j  \right) \\
U_1 &=& U_0 + \Delta t ( b_1 K_1 + \ldots + b_s K_s)
\end{array}
\end{equation}
has order $p$. Compared to Singly-RKTASE methods, instead of
$s$ coefficients $\beta_i$ we have $s r$ coefficients $\beta_{ij}$, providing more freedom to improve accuracy and stability properties. Although we could take different $r_i$ for each stage, we have taken the same $r$ for simplicity. Note that if we impose that all operators $T_i$ in \eqref{A10} have order $p$ with $r=p$, the coefficients $\beta_{ij}$ are uniquely determined, resulting in a Singly-RKTASE. Then, we must consider $r>p$, which increases computational cost. Alternatively, it is possible to use operators $T_i$ with order $q < p$ in such a way that the resulting modified Singly-RKTASE method has order $p$. This requires further analysis of the conditions the methods's  coefficients must satisfy to achieve order $p$.

The rest of the paper is organized as follows: In section 2, we prove that modified Singly-RKTASE methods are  a particular case of W-methods. In section 3, we derive the conditions that the coefficients $\beta_{ij}$ and $\alpha$ of the MSRKTASE must satisfy to have order $p$. In section 4 we develop MSRKTASE methods of orders 2 and 3 with optimal stability and accuracy properties. Finally, in section 5,
through numerical experiments, we demonstrate the effectiveness of the newly developed methods, highlighting their potential for solving both linear and nonlinear stiff systems with improved performance.

\section{Equivalence between modified Singly-RKTASE and W-methods}

Given a matrix $W$ (usually some approximation to the Jacobian matrix at the current integration point $(t_n, Y_n)$) a singly $W$-method \cite{SteWolf:79} with $\hat s$ stages provides an approximation to the solution of Eq.~(\ref{A1}) at  $t_{n+1}=t_n + h$   by the formulas
\[
\begin{array}{l}
(I-h\alpha W) \widehat K_1 =h f(Y_n) \\
(I-h\alpha W) \widehat K_i = h f( Y_n+ \sum_{j=1}^{i-1} \hat a_{ij} \widehat K_j)+
    h\alpha W \sum_{j=1}^{i-1} l_{ij} \widehat K_j, \qquad i=2, \ldots, \hat s \\
    Y_{n+1} = Y_n + \sum_{i=1}^{\hat s} \hat b_i \widehat K_i
    \end{array}
\]

The coefficients of the method can be arranged in a Butcher tableau as
\[
\begin{array}{c|c|c}
\hat c & \hat A & L \\
\hline
 & \hat b^T & \alpha
 \end{array}\quad = \quad
\begin{array}{c|cccc| ccccc}
0 &  & & & & &
\\
\hat c_2  & \hat a_{21} & & & &  l_{21}
\\
\vdots & & \ddots& & & \vdots &  \ddots
\\
\hat c_s  & \hat a_{s1} & \ldots & \hat a_{s,s-1} & &
l_{s1} & \ldots &l_{s,s-1}
\\
\hline
 & \hat b_1 & \ldots & \hat b_{s-1} & \hat b_s & \alpha
 \end{array}
\]

On the other side, the product of the operator $T_i$ in \eqref{A10} by a vector $v$ can be written as
\[
\begin{array}{l}
w_1 = (I-h\alpha W)^{-1} v \\
w_j =(I-h\alpha W) w_{j-1}, \quad j=2, \ldots, r \\
T_i v = \beta_{i1} w_1 + \ldots + \beta_{ir} w_r
\end{array}
\]
Then, the equations \eqref{A10}\eqref{A11} defining a Modified Singly-RKTASE scheme can be written as
\[
\begin{array}{l}
(I-h\alpha W) \widehat K_1 = h f(Y_n) \\
(I-h\alpha W) \widehat K_j = K_{j-1}, \quad j=2, \ldots, r \\
K_1= \beta_{11} \widehat K_1 + \ldots + \beta_{1r} \widehat K_r \\
(I-h\alpha W) \widehat K_{r(i-1)+1} = h f( Y_n+ \sum_{j=1}^{i-1} a_{ij} K_j)  \qquad    i= 2, \ldots,s\\
(I-h\alpha W) \widehat K_{r(i-1)+j} = K_{r(i-1)+j-1}, \quad j=2, \ldots, r \\
K_i= \beta_{i1} \widehat K_{r(i-1)+1} + \ldots + \beta_{ir} \widehat K_{r i } \\
Y_{n+1}= Y_n + \sum_{i=1}^s b_i K_i
\end{array}
\]

Taking into account that
$
\widehat K_j = \widehat K_{j-1} + h \alpha W \widehat K_j,
$
we have
\[
\begin{array}{l}
(I-h\alpha W) \widehat K_1 = h f(Y_n) \\
(I-h\alpha W) \widehat K_j = h f(Y_n) + h\alpha W \sum_{l=1}^{j-1} K_{l}, \quad j=2, \ldots, r \\
K_1= \beta_{11} \widehat K_1 + \ldots + \beta_{1r} \widehat K_r \\
(I-h\alpha W) \widehat K_{r(i-1)+1} = h f( Y_n+ \sum_{j=1}^{i-1} a_{ij} K_j)  \qquad    i= 2, \ldots,s\\
(I-h\alpha W) \widehat K_{r(i-1)+j} = h f( Y_n+ \sum_{j=1}^{i-1} a_{ij} K_j) + h\alpha W \sum_{l=1}^{j-1} K_{r(i-1)+l-1}, \quad j=2, \ldots, r \\
K_i= \beta_{i1} \widehat K_{r(i-1)+1} + \ldots + \beta_{ir} \widehat K_{r i } \\
Y_{n+1}= Y_n + \sum_{i=1}^s b_i K_i
\end{array}
\]
and this leads to
\[
\begin{array}{l}
(I-h\alpha W) \widehat K_1 = h f(Y_n) \\
(I-h\alpha W) \widehat K_j = h f(Y_n) + h\alpha W \sum_{l=1}^{j-1} K_{l}, \quad j=2, \ldots, r \\
(I-h\alpha W) \widehat K_{r(i-1)+j} = h f( Y_n+ \sum_{j=1}^{i-1} a_{ij} \sum_{l=1}^r \beta_{jl} \widehat K_{r(i-1)+j}) + h\alpha W \sum_{l=1}^{j-1} K_{r(i-1)+l-1}, \quad  i=2, \ldots, s, \quad j=2, \ldots, r \\
Y_{n+1}= Y_n + \sum_{i=1}^s b_i \sum_{l=1}^r \beta_{il} \widehat K_{r(i-1)+j}
\end{array}
\]

This is just an $W$ method with $sr$ stages whose vector $\widehat b$ and matrices $\widehat A$, $L$ are given by
\[
L= \begin{pmatrix}
L_r & & \\
 & \ddots \\
 &  & L_r
\end{pmatrix}\in \mathbb{R}^{rs\times rs},
\qquad
\hat A =
\begin{pmatrix}
0 & & \\
A_{21} & 0\\
\vdots & \\
A_{s1} & \cdots & A_{s,s-1} & 0
\end{pmatrix}\in \mathbb{R}^{rs\times rs},
\]
\[
\widehat b = (b_1 \beta_{11},\ldots, b_1 \beta_{1r}, \ldots, b_s \beta_{s1},\ldots, b_s \beta_{sr})^T
\]
with
\[
L_r=
\begin{pmatrix}
0 & &  \\
1& 0 & \\
\vdots & \ddots & \ddots\\
1 & \cdots & 1 & 0
\end{pmatrix}\in \mathbb{R}^{r\times r}, \qquad
A_{i,j} = a_{ij}
\begin{pmatrix}
 \beta_{j1} &  \ldots &  \beta_{jr} \\
\vdots & \\
\beta_{j1} & \ldots & \beta_{jr}
\end{pmatrix}\in \mathbb{R}^{r\times r}
\]
This equivalence let us analyze the order of the Modified Singly-RKTASE
through the order conditions of W-methods. Also, the absolute stability can be studied through the $W$-methods.
The stability function of an $W$ method is given by
\[
\hat R(z)= 1 + z \hat b^T (I-z(\hat A + \Gamma))^{-1} \mathbf{1}
\]
and the limit of this function when $z$ goes to infinity is given by
\begin{equation}
\label{rinfinity}
\lim_{z\to\infty} \hat R(z) = 1-\hat b^T(\hat A + \Gamma)^{-1} \mathbf{1}.
\end{equation}
We will use this to develop Modified Singly-RKTASE methods. Recall that
\begin{itemize}
\item A method is called A-stable if $\vert \hat R( z)\vert\le 1$ for all Re $z\le 0$;
\item
A method is said to be L--stable if it is A--stable and
 $\hat R(\infty)=0$;
\item
A method is called A$(\theta)$--stable if $\vert \hat R(z)\vert<1$
for all $z$ such that $\arg(-z) \le \theta$, that is, its stability
region contains the left hand region of the complex plane with angle $\theta$.
If in addition $\hat R(\infty)=0$,
it is called L$(\theta)$--stable.
\end{itemize}

\section{Order conditions}
Denoting $\Gamma = \alpha(I+L)$, and $\mathbf{1}=(1,\ldots, 1)^T$,
the order conditions for an W-method are given in Table \ref{table1}
(see e.g. \cite{HaWa:96, SteWolf:79}).
\begin{table}
\caption{order conditions for $W$-methods up to order 4}
\label{table1}
\[
\begin{array}{|ll|lll|}
\hline
\hat b^T \mathbf{1} =1& & && \\
\hline
\hat b^T \hat c= 1/2  &  & \hat b^T \Gamma \mathbf{1} =0 & & \\
\hline
\hat b^T \hat c^2= 1/3 & \hat b^T\hat A \hat c=1/6  &  \hat b^T \Gamma^2 \mathbf{1} =0 &  \hat b^T \hat A \Gamma \mathbf{1} =0 &  \hat b^T \Gamma \hat A \mathbf{1} =0\\
\hline
\hat b^T \hat c^3= 1/4 & \hat b^T\hat A \hat c \cdot \hat c=1/8  &  \hat b^T \Gamma^3 \mathbf{1} =0 &  \hat b^T \hat A \Gamma^2 \mathbf{1} =0 &  \hat b^T \Gamma \hat A \Gamma \mathbf{1} =0\\
 \hat b^T\hat A \hat c^2=1/12  &  \hat b^T\hat A^2 \hat c =1/24 &
 b^T \Gamma^2 \hat A \mathbf{1} =0 &  \hat b^T \hat A^2 \Gamma \mathbf{1} =0 &  \hat b^T \hat A \Gamma \hat A \mathbf{1} =0 \\
 & & \hat b^T \Gamma \hat A^2 \mathbf{1} =0 & \hat b^T \Gamma \hat c^2 =0 & \hat b^T(\hat A\Gamma\mathbf{1}) \cdot \hat c=0\\
 \hline
\end{array}
\]
\end{table}
In the case the operator $W$ is exactly the Jacobian matrix, $W=\partial f/\partial y (t_n, Y_n)$ (Rosenbrock-Wanner methods), many elementary differentials
become the same and the order conditions reduce to those in Table \ref{table2}.

\begin{table}
\caption{Order conditions for a Rosenbrock method ($W=\partial f/\partial y$) up to order 4}
\label{table2}
\[
\begin{array}{|l|lll|}
\hline
\hat b^T \mathbf{1} =1& & & \\
\hline
\hat b^T (\Gamma+\hat A) c =1/2 & & &\\
\hline
\hat b^T \hat c^2= 1/3 & \hat b^T(\Gamma+\hat A)^2 \mathbf{1}=1/6  & & \\
\hline
\hat b^T \hat c^3= 1/4 & \hat b^T(\hat A(\Gamma+ \hat A)\mathbf{1}) \cdot \hat c=1/8  &
 \hat b^T(\Gamma+\hat A) \hat c^2=1/12 &
   \hat b^T ( \Gamma +\hat A) ^3 \mathbf{1} =1/24 \\
 \hline
\end{array}
\]
\end{table}

Let us analyze the order conditions for a Modified Singly-RKTASE scheme expressed as its equivalent $W$-method.
We can take without losing generality $\beta_{i1}+ \ldots + \beta_{ir}=1$ because this is equivalent to re-scale the coefficients $b, A$ of the underlying RK method.

First, the order conditions that do not involve the matrix $\Gamma$ reduce to the order conditions associated to the undelying RK method. For example, $\hat b^T \mathbf{1}=b_1+ \ldots+b_s$, or $\hat b^T \hat A \hat c = b^T A c$ which are the corresponding equations of the RK method.

On the other hand, the form of the matrix $L$ introduces some
incompatibility between some of the order conditions leading to the next
\begin{theorem}
A Modified Singly-RKTASE method can not have order $r+1$.
\end{theorem}
\Proof
Since the matrix $\Gamma$ is a block diagonal matrix with blocks $\alpha(I_r+L_r)$, it is clear that $(\Gamma- \alpha I)^r =(\alpha L)^r=0$ because
$L_r$ is strictly lower triangular with dimension $r$. Then,
$(\Gamma- \alpha I)^r \mathbf{1}=0$. If the method had order $r+1$, by the order conditions it should be $\widehat b^T(\Gamma- \alpha I)^r \mathbf{1}=\alpha^r$ which is not zero if $\alpha \ne 0$.
\Endproof

With this result, a Modified Singly-RKTASE method with order $p$ must have
necessarily $r \ge p$.

\section{Modified Singly-RKTASE methods with orders 2 and 3 }

\subsection{Methods with $s=2$, $r=2$ and order 2}
 There exists a family of Runge-Kutta
schemes with two stages and order two that depend on one coefficient $c_2$ given by
\[
a_{21}=c_2,\quad
b_2=\dfrac{2 - 3 c3}{6 c2 (c2 - c3)}, \quad
b_1=1-b_2
\]
Since $b^T A c=0$ for any $c_2$, the associated order condition of order three does not depend on $c_2$ and the other one vanishes for $c_2= 2/3$. So, we will take this value of $c_2$.

The only remaining condition for a Modified Singly-RKTASE  with $s=r=2$ to achieve order 2 is
$\hat b^T \Gamma \mathbf{1} =0$
which is satisfied for
\[
\beta_{22}= -(4+\beta_{12})/3.
\]
The limit of the stability function when $z$ approaches infinity for the case of our Modified Singly-RKTASE method of order two is given by
\[
\lim_{z\to\infty} \hat R(z) = 1-\hat b^T(\hat A + \Gamma)^{-1} \mathbf{1} =
-\dfrac{-7 + 12 \alpha - 6 \alpha^2 + 6 \beta_{12} + \beta_{12}^2}
 {6 \alpha^2}.
\]
There are two values of $\beta_{12}$ that make $\hat R(\infty)=0$. The one that yields better results in terms of stability and accuracy is
\[
\beta_{12}= -3 + \sqrt{16 - 12 \alpha + 6 \alpha^2}
\]
There is only one free parameter $\alpha$ left, which we will use to maximize the stability region and minimize the error coefficients. The norm of the error coefficients of order 3 is given by
\[
C_3=\left(
(\hat b^T\hat A \hat c-1/6)^2+
(\hat b^T \Gamma^2 \mathbf{1})^2+ (\hat b^T \hat A \Gamma \mathbf{1})^2+ (\hat b^T \Gamma \hat A \mathbf{1})^2 \right)^{1/2}.
\]
The corresponding norm when $W= f'$ is
\[
D_3= \left|\hat b^T(\Gamma+\hat A)^2 \mathbf{1}-1/6 \right|.
\]
The method is $A$-stable
(and therefore $L$-stable) if $\alpha \in [ 0.3117, 3.257]$ and $C_3$
is monotonic decreasing with $\alpha$. A good value of the parameter can is $\alpha=32/100$. With these values, we obtain an $L$-stable method of order 2, which we will name MSRKTASE2.  The error coefficients are given in
Table \ref{tabla1}. For comparison, we also include the error coefficients of the Singly-RKTASE
of order 2 obtained in \cite{CaFuMoRa:23}, named SRKTASE2. As we can see, the new method achieves $L$-stability
(SRKTASE2 has $\hat R(\infty)=1/2$) and error coefficients that are about 20 times smaller (40 times smaller in the case where $W=f'$).

We could use the parameter $\beta_{12}$ to obtain smaller error coefficients without imposing the condition $R(\infty)=0$. For example we could nullify the coefficient  $D_3$
(order three for the case $W=f'$,) but in this case,
we only achieve reasonable stability regions for large values of $\alpha$ and $C_3$ becomes very large (greater than one hundred).

\begin{table}
\caption{properties of the proposed methods}
\label{tabla1}
\begin{center}
 \begin{tabular}{|l|l|l|l|l|l|}
\hline
Method & $p$ & $\alpha$ & $C_{p+1}$  & $D_{p+1}$ ( $W= f'$ ) & $\theta$ \vrule height 12pt width 0pt\\[2pt]
\hline
SRKTASE2 & $2$ & $2$ & $4.00347 $  &  $4.16667$ & $90^{\circ} $  \vrule height 12pt width 0pt\\[2pt]
\hline
MSRKTASE2 & $2$  & $0.32$ &  $0.212866$ &  $0.10116$ &  $90^{\circ} $ \vrule height 12pt width 0pt\\[2pt]
\hline
SRKTASE3 & $3$ & $1.8868$ & $6.7171 $  &  $6.6753$ & $88.99^{\circ} $  \vrule height 12pt width 0pt\\[2pt]
\hline
MSRKTASE3a & $3$  & $0.54$ &  $0.1817$ &  $0.2288$ &  $88.23^{\circ} $ \vrule height 12pt width 0pt\\[2pt]
\hline
MSRKTASE3b & $3$ & $0.56 $ &  $0.3968$  &  $0.0035 $  &  $50.38 ^{\circ} $\vrule height 12pt width 0pt\\[2pt]
\hline
\end{tabular}
\end{center}
\end{table}

\subsection{Methods with $s=3$, $r=3$ and order 3}
 There exist a family of Runge-Kutta
schemes with three stages and order three, defined by the coefficients
\[
a_{32}=\dfrac{c_3(c_2-c_3)}{c_2(3c_2-2)},\quad
b_2=\dfrac{2 - 3 c3}{6 c2 (c2 - c3)}, \quad
b3=\dfrac{2 - 3 c2}{6 c3 (c3 - c2)}, \quad
b_1=1-b_2-b_3
\]
depending on two free parameters $c_2$ and $c_3$ with $c_2\ne c_3$ and
$c_2 \ne 2/3$.

The remaining conditions for a Modified Singly-RKTASE  with $s=r=3$ to achieve order 3 form a set of four equations
$\hat b^T \Gamma \mathbf{1} =0$,
$\hat b^T \Gamma^2 \mathbf{1} =0$, $\hat b^T \hat A \Gamma \mathbf{1} =0 $,  $\hat b^T \Gamma \hat A \mathbf{1} =0$
which can be solved for the parameters $\beta_{12}, \beta_{13}, \beta_{23}$ and $\beta_{33}$, resulting in:
\[
\begin{array}{l}
\beta_{12}=\dfrac{c_3 (-2 + 3 c_3) \beta_{22} - 3 c_2^2 (6 c_3 + \beta_{32}) +
  2 c_2 (9 c_3^2 + \beta_{32})}{(c_2 - c_3) (2 - 3 c_3 +
    c_2 (-3 + 6 c_3))},  \\[9pt]
\beta_{13}= -\dfrac{c_3 (-2 + 3 c_3) (1 + \beta_{22}) -
      3 c_2^2 (1 + 4 c_3 + \beta_{32}) +
      2 c_2 (1 + 6 c_3^2 + \beta_{32})}{2 (c_2 - c_3) (2 - 3 c_3 +
        c_2 (-3 + 6 c_3))}\\[7pt]
\beta_{23} = 1/2 (-1 - \beta_{22}), \\
\beta_{33} = 1/2 (-1 - \beta_{32}),
\end{array}
\]
depending on the parameters $c_2, c_3, \beta_{22}$ and $\beta_{32}$. Note that in this family of methods of order three, $\alpha$ is also a free parameter.
Imposing $\beta_{32}=\beta_{22}=\beta_{12}$ we obtain the third order Singli-RKTASE methods obtained in \cite{CaFuMoRa:23}.

The free parameters can be selected to maximize the stability region and to minimize the coefficients of the leading term of the local error.
For these methods of order three, it is also satisfied that
\[
\hat b^T\hat A \Gamma \hat A \mathbf{1}=0, \quad
\hat b^T\hat A^2  \Gamma\mathbf{1}=0, \quad
\hat b^T \Gamma \hat A^2  \mathbf{1}=0,\quad
\hat b^T \Gamma \hat A \Gamma \mathbf{1}=0, \quad
\hat b^T \Gamma \hat c^2=0,\quad \hat b^T(\hat A\Gamma\mathbf{1}) \cdot \hat c=0
\]
and the other error coefficients of the term of order four are
\begin{equation}
\label{eqrest}
\begin{array}{l}
C_{41}\equiv\hat b^T\hat A^2 \hat c -1/24 =  b^T A^2  c -1/24 = -1/24 \\[5pt]
C_{42}\equiv(8\hat b^T\hat A \hat c \cdot \hat c -1)/24 = (8 b^T A  c \cdot c-1)/24 =(4c_3-3)/72 \\[5pt]
C_{43}\equiv(12 \hat b^T\hat A \hat c^2 -1)/24 = (12 b^T A  c^2 -1)/24 =(2c_2-1)/24 \\[5pt]
C_{44}\equiv(4\hat b^T\hat c^3-1)/24  =  (4b^T c^3-1)/24  = (2((c_2 (2 - 3 c_3) + 2 c_3)-3)/72 \\[6pt]
C_{45}\equiv\hat b^T\hat A \Gamma^2 \mathbf{1}=
\alpha^2 \frac{-2 \beta_{22} + 6 c_3 (3 + \beta_{22}) - 3 c_3^2 (3 + \beta_{22}) + 2 \beta_{32} +
     9 c_2^2 (3 + \beta_{32}) -
     3 c_2 (6 - \beta_{22} + 2 c_3 (3 + \beta_{22}) + 3 \beta_{32}))}{12 (c_2 - c_3) (2 - 3 c_3 +
     c_2 (-3 + 6 c_3))} \\[6pt]
C_{46}\equiv\hat b^T \Gamma^2 \hat c= \alpha^2 \frac{(-2 \beta_{22} + 3 c_3 (3 + \beta_{22}) + 2 \beta_{32} - 3 c_2 (3 + \beta_{32}))}{12 (c_2 - c_3)} \\[7pt]
C_{47}\equiv\hat b^T \Gamma^3 \mathbf{1}= \alpha^3
\end{array}
\end{equation}

If we solve $\hat b^T \Gamma^2 \hat c=0$ and $\hat b^T\hat A \Gamma^2 \mathbf{1}=0$ for $\beta_{22}$ and $\beta_{32}$ we obtain the method in \cite{CaFuMoRa:23}, independent of the values of $c_2$ and $c_3$.

The first four equations in \eqref{eqrest} depend only on $c_2$ and $c_3$
and can not be satisfied simultaneously. The 2-norm of these four equations reaches its minimum value at $c_2= 0.496188$, $c_3=0.764887$, very close to $c_2=1/2$, $c_3=3/4$ for which
\[
\dfrac{\big((24 b^T A^2c-1)^2+(12b^T Ac^2-1)^2+(8b^T Ac\cdot c-1)^2+(4b^T c^3-1)^2\big)^{1/2}}{24}=
\dfrac{\sqrt{145}}{288}=0.041811
\]
From now on, we will take $c_2=1/2$ and $c_3=3/4$. Thus, we have three coefficients $\alpha, \beta_{22}, \beta_{32}$ to minimize the error coefficients and maximize the stability region.

The limit of the stability function of this method as $z$ approaches infinity is given by
\[
\lim_{z\to\infty} \hat R(z) = 1-\hat b^T(\hat A + \Gamma)^{-1} \mathbf{1} =
\dfrac{a_0 + a_1\alpha -288 \alpha^2 +96\alpha^3}{96\alpha^3}
\]
with
\[
\begin{array}{l}
a_0= -(-3 + \beta_{22}) (-3 + \beta_{32}) (33 + 3 \beta_{22} + 4 \beta_{32})\\
a_1= -6 (-45 + 12 \beta_{22} + \beta_{22}^2)
\end{array}
\]
There are two values of $\beta_{32}$ (as functions of $\alpha$ and $\beta_{22}$) that make
$R(\infty)=0$. We will select one of these values and use the coefficients $\alpha$ and $\beta_{22}$ to minimize the error coefficients and maximize the $L(\theta)$-stability angle.

The 2-norm of the error coefficients is given by
\[
C_4=( C_{41}^2+ \cdots + C_{47}^2)^{1/2}
\]
Minimizing this norm is equivalent to minimizing (the other coefficients are constant)
$
C_{45}^2 + C_{46}^2 + C_{47}^2.
$
A good compromise between large $\theta$ and small $C_4$ is obtained
for
\[
\alpha = 0.54, \quad \beta_{22} = -6.1,\quad  \beta_{32} = -2.75034
\]
resulting in a method, that we will name MSRKTASE3a,  with stability angle and error coefficients given in
Table \ref{tabla1}.

Alternatively, we can minimize the 2-norm of the coefficients of the leading term in the case that  $W$ is exactly the Jacobian matrix,
\[
D_4=( D_{41}^2+ \cdots + D_{45}^2)^{1/2}.
\]
This is equivalent to minimize $|D_{41}|= \hat b^T ( \Gamma +\hat A) ^3 \mathbf{1} -1/24$ (the other coefficients are constant). This error coefficient can be made to vanish when
\[
\beta_{22}=-3 - \dfrac{1}{3 \alpha^2} + 8 \alpha.
\]
A good compromise between maximizing $\theta$ and minimizing $C_4$ with $\alpha$ is obtained when
\[
\alpha = 0.56, \quad \beta_{22} = -6.1, \quad \beta_{32} = -2.75034
\]
resulting in a method, that we will name MSRKTASE3b,  with stability angle and error coefficients given in
Table \ref{tabla1}.

\begin{figure}
\includegraphics[width=0.49\textwidth]{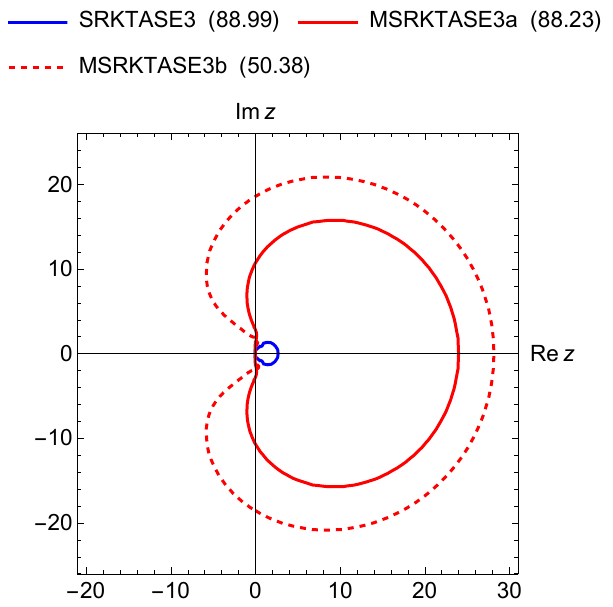}
\  \  \
\includegraphics[width=0.48\textwidth]{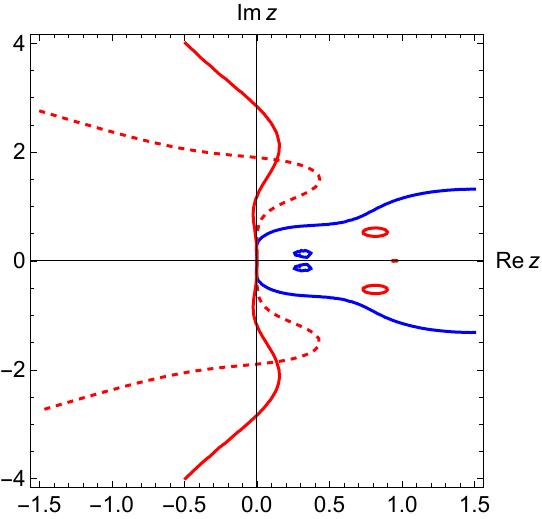}
\caption{Stability regions of the new methods (left). On the right,
a zoom of the region near the origin.}
\label{figura1}
\end{figure}

In Figure \ref{figura1} we plot the boundary of the stability regions of the proposed methods. The method with minimal
$D_4$ is shown in red-dashed, the method with minimal $C_4$ in red-solid, and the method from \cite{CaFuMoRa:23} in blue.

We can see from these results that the new method
MSRKTASE3a, has similar stability properties to the method of order 3 in \cite{CaFuMoRa:23}, but has about 35 times smaller error coefficients. The new method is expected to provide more accurate approximations with the new method.  The other method, MSRKTASE3b, also has much smaller
error coefficients, although the stability angle is smaller. Nonetheless, its stability region is not significantly worse.

\section{Numerical experiments}

To evaluate the performance of the new methods, we considered two test problems and integrated them using the two new methods, MSRKTASE3a and MSTASE3b. For comparison, we also used the method proposed in \cite{CaFuMoRa:23} (STASE) to demonstrate the improvements of these new modified singly-RKTASE methods over the previous singly-RKTASE methods. Additionally, to benchmark the performance of the new methods against other known Runge-Kutta methods for stiff problems, we integrated the problems with a Singly-Diagonally Implicit RK method of order 3 (SDIRK) proposed in \cite{KeCa:16}. Since the SDIRK method is implicit,
we used a simplified Newton method
to solve the stage equations, approximating the Jacobian matrix with the matrix $W$ used in the TASE operators.

For each method and problem, we computed the CPU time required for solving it and the 2-norm of the global error at the end of the integration interval. These data points allowed us to plot the global error (in logarithmic scale) against the CPU time, resulting in an efficiency plot.

\noindent {\bf Problem 1:} (taken from \cite{Bassenne:2020} and \cite{CaMoRa:21}) The 1D diffusion of a scalar function $y = y(x, t)$, with a time dependent source term
\[
\dfrac{\partial y}{\partial t}= \dfrac{\partial^2 y}{\partial x^2}
+ 0.1 \sin (t/50), \qquad y(x, 0) = 1 - \cos(x)^{101}, \qquad  0 \le x \le 2\pi.
\]
The solution $y = y(x, t)$ is assumed to be $2\pi$-periodic in $x$.
For the spatial discretization, we used fourth-order centered difference schemes
with a grid resolution of $N = 512$, assuming periodic boundary conditions. The real part of the eigenvalues of the Jacobian matrix of the semi-discrete problem ranges from $-3.5 \times 10^4$ to $2.79 \times 10^{-5}$.

\begin{figure}
\begin{center}
\includegraphics[width=0.5\textwidth]{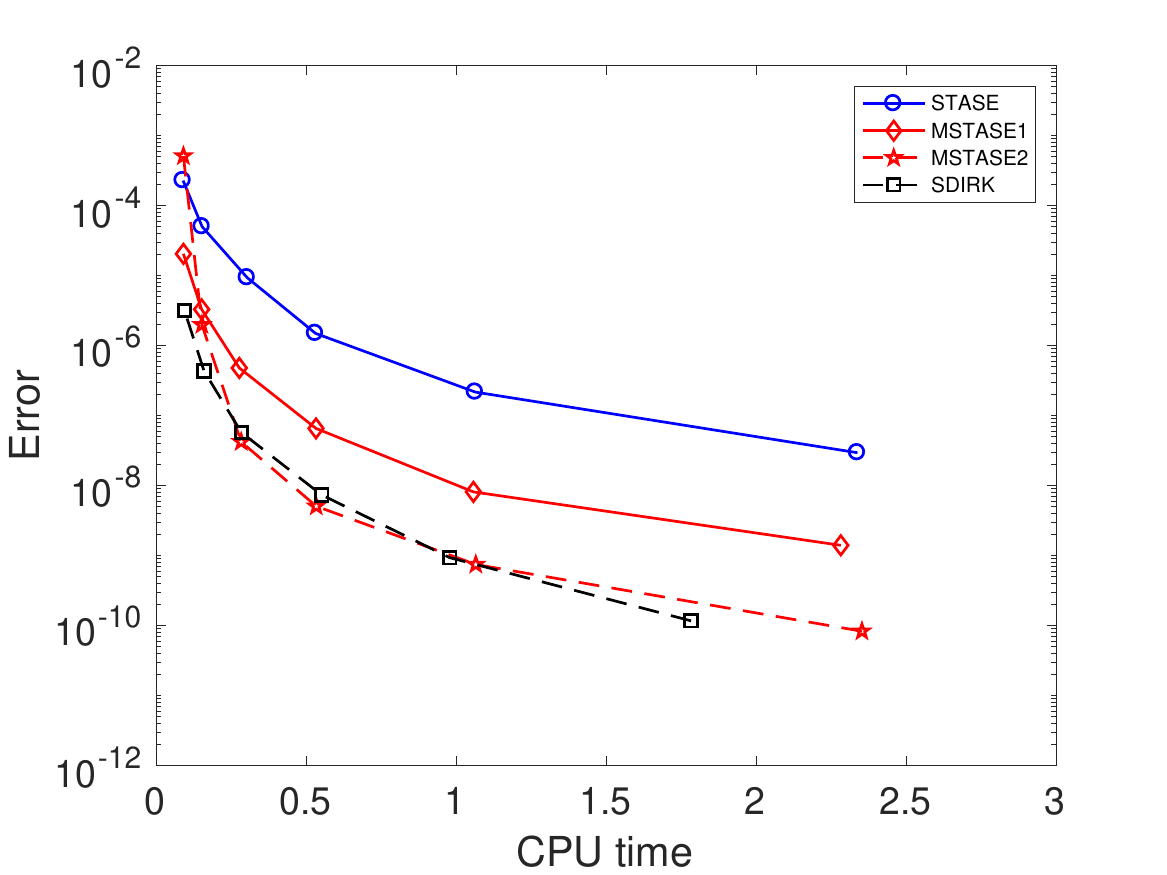}
\caption{Efficiency plot for Diffusion problem}
\label{problem1}
\end{center}
\end{figure}

The efficiency plot of the methods for this problem is depicted in
Figure \ref{problem1}.
The method MSRKTASE3b is as efficient as SDIRK and more efficient than
MSRKTASE3a. This is due to the fact that the eigenvalues of the Jacobian matrix are close to the real axis, therefore included in the stability region and the Jacobian matrix $J$ of the problem is constant, then $W=J$.  In this situation MSRKTASE3b has smaller error coefficients than MSRKTASE3a and this translates into a better efficiency. MSRKTASE3a is clearly more efficient than SRKTASE.

\bigskip

\noindent {\bf Problem 2:} The 1D Burgers’ equation written in the conservative
form
\[
\dfrac{\partial y}{\partial t}=0.1 \dfrac{\partial^2 y}{\partial x^2}
-\dfrac{\partial}{\partial x} \left(\dfrac{y}{2}\right)^2, \qquad  y(x,0)= 1 -\cos(x)^{101}, \qquad  0 \le x \le 2\pi.
\]
The solution $y=y(x,t)$ is assumed to be $2\pi$-periodic in $x$.
For the spatial discretization,  a fourth-order centered difference scheme with grid resolution of $N=512$.

The real part of the eigenvalues of the Jacobian matrix of the semi-discrete problem ranges from $-3.5\times 10^{3}$ to $2.6\times 10^{-6}$.

\bigskip

We integrated Burger's problem using the three methods with three different options for the matrix $W$.   Firstly, the Jacobian matrix is evaluated at every time step and $W= \partial f(t_n, y_n)/\partial y $. With this option, the $LU$ matrix factorization had to be computed at every step, increasing the computational cost. Secondly, the Jacobian matrix was evaluated only at the initial time step and $W= \partial f(t_0, y_0)/\partial y $, requiring only one $LU$ factorization and considerably reducing the CPU time. However, this could affect the accuracy in the case of TASE methods or the number of iterations required to solve the non linear system in the case of the SDIRK method. Finally, $W$ was taken as the linear part of the semi-discrete differential equation (the matrix of the diffusion term). The computational cost is similar to the second option, but the accuracy could be lower. This option provided insight into the methods behaviour when $W$ poorly
approximates the Jacobian matrix.

\begin{figure}
\begin{center}
\includegraphics[width=0.495\textwidth]{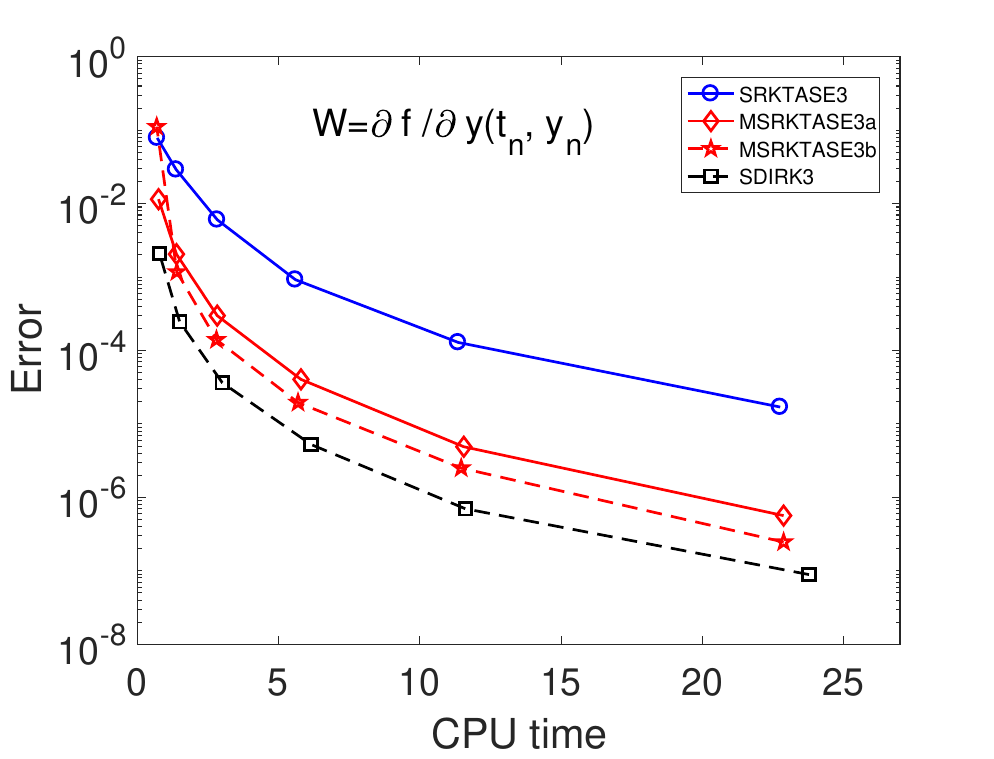}\
\includegraphics[width=0.495\textwidth]{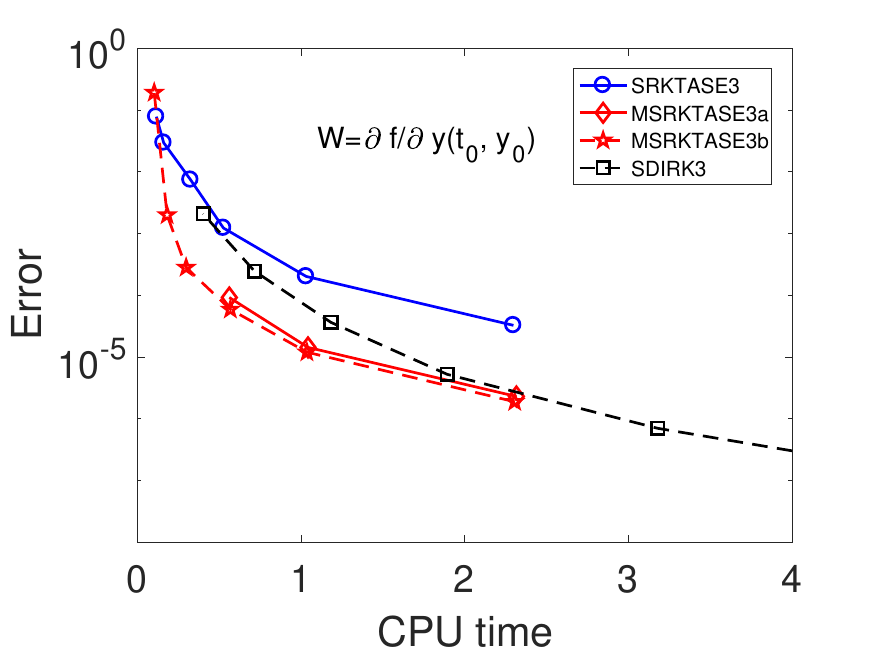}
\end{center}
\caption{Efficiency plot for Burger's problem integrated
with $W$ the Jacobian matrix evaluated at every step (left) and
$W$ the Jacobian matrix evaluated just at the initial time step (right)}
\label{problem2a}
\end{figure}

\begin{figure}
\begin{center}
\includegraphics[width=0.5\textwidth]{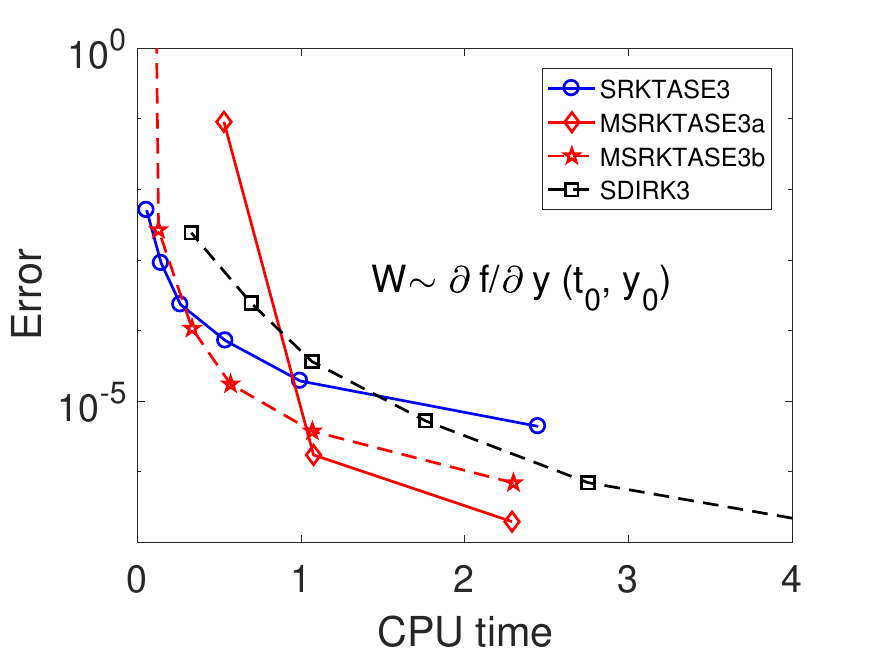}
\end{center}
\caption{Efficiency plot for Burger's problem integrated
with $W$ equal to the matrix of diffusion term}
\label{problem2b}
\end{figure}

Figures \ref{problem2a} and \ref{problem2b} show the efficiency plots for Burger's problem integrated with these three options.

As shown in the Figures, evaluating the Jacobian matrix only at the initial step reduces the CPU time but increases the global error for the TASE-based methods, where the local error depends on $W$ and is smaller when it coincides with the Jacobian matrix.  For the SDIRK method, the error does no vary significantly with different $W$, but solving the nonlinear system requires more iterations when the Jacobian matrix is less accurately approximated, thereby increasing the computational cost.

MSRKTASE3a exhibited stability issues with the two largest time steps when $W$ is not the Jacobian matrix at every step, likely due to a high sensitivity of the stability function to parameter changes. MSRKTASE3b also showed stability problems with the largest time step size and the poorest approximation of the Jacobian matrix. Further research is being carried out in this regard.

Concerning the efficiency of the methods, when the jacobian matrix is evaluated at every step, SDIRK3 is the most efficient, followed by MSRKTASE3b and MSRKTASE3a. SRKTASE3 is the least efficient due to its larger error coefficients. When the Jacobian matrix is evaluated only at the initial point, MSRKTASE3a and MSRKTASE3b are more efficient while SRKTASE3 remains the least  efficient.  Finally, when the matrix $W$ is the matrix of the diffusion term, MSRKTASE3a is more efficient for small time steps, where it has no stability problems, due to the relevance of the error coefficient C4. MSRKTASE3b remains more efficient than SDIRK and SRKTASE3 is also more efficient than SDIRK except for the smallest time step. The solution of the nonlinear equations makes SDIRK less efficient.

\section{Conclusions}
We presented a modification of the Singly-RKTASE methods for the numerical solution of stiff differential equations by taking different TASE operators for each stage of the RK method. We proved that these methods are equivalent to $W$-methods which enabled us to derive the order conditions when the TASE operators have order smaller than the order $p$ of the Runge-Kutta scheme. For the case where $p=3$ we obtained Modified Singly-RKTASE methods that significantly reduce the error coefficients compared to Singly-RKTASE methods, while maintaining stability.

Numerical experiments on both linear and nonlinear stiff systems demonstrate that the modified Singly-RKTASE methods provide significant improvements in accuracy and computational efficiency over the original Singly-RKTASE schemes, making them competitive with other methods such as Diagonally Implicit RK methods.

{\small

\end{document}